\documentclass[11pt]{article}
\usepackage{amsmath,amsfonts,amssymb,hyperref,setspace}
\newtheorem{thm}{Theorem}[section]
\newtheorem{lem}[thm]{Lemma}
\newtheorem{cor}[thm]{Corollary}
\newtheorem{exa}[thm]{Example}
\newtheorem{defn}[thm]{Definition}

\numberwithin{equation}{section}
\title{ Metrizability of Cone Metric spaces}
\author{{}\\{\normalsize{\sc
${}$\, Mehdi Asadi$^{a,}$\thanks{Corresponding author.\quad Fax:+98-241-4220030.},
${}$\, S. Mansour Vaezpour$^{b}$,
${}$\, Hossein Soleimani$^{c}$ }}\\
   \\{\footnotesize{\it${}^{a}$  Islamic Azad University, Zanjan Branch, Zanjan, Iran}}
  \\{\footnotesize{\it { masadi.azu@gmail.com}}}
  \\{\footnotesize{\it${}^b$ Dept. of Math., Amirkabir University of Technology, Tehran, Iran}}
     \\{\footnotesize{\it {vaez@aut.ac.ir}}}
    \\{\footnotesize{\it${}^{c}$  Islamic Azad University, Malayer Branch, Malayer, Iran}}
    \\{\footnotesize{\it {hsoleimani54@gmail.com}}}
        }
        \date{\empty}
\begin{document}
\maketitle
%---------------------------------------------------------------------------------------------------------------------
                                                        %abstract
%---------------------------------------------------------------------------------------------------------------------

\begin{abstract}
In  2007 H. Long-Guang and Z. Xian, [H. Long-Guang and Z. Xian, \emph{ Cone Metric Spaces and Fixed Point Theorems of Contractive Mapping}, J. Math. Anal. Appl.,
\textbf{322}(2007), 1468-1476],
 generalized the concept of a metric space, by introducing cone metric spaces,
 and obtained some fixed point theorems for mappings satisfying certain contractive conditions.
 The main question was "Are cone metric spaces a real generalization of metric spaces?" Throughout this paper we answer the question in the negative, proving that every cone metric space is metrizable  and the equivalent metric satisfies the same contractive conditions as the cone metric. So most of the fixed point theorems which have been proved are straightforward results from the metric case.

\end{abstract}

\section{Introduction and Preliminary}
 Long-Guang and Xian in \cite{Xian} generalized the concept of a metric space, replacing the set of real numbers by an ordered Banach space and obtained some fixed point theorems for mapping satisfy different contractive conditions.\\

Recently Wei-Shih Du in \cite{Du1} has proved that the Banach contraction principle in general metric spaces and in TVS-cone metric spaces are equivalent, and in \cite{Du2} has obtained new type fixed point theorems for nonlinear multivalued maps in metric spaces and the generalizations of Mizoguchi-Takahashi's fixed point theorem and Berinde-Berinde's fixed point theorem. But in this paper we introduce a equivalent metric which satisfies the same contractive conditions as the cone metric in the easy way.

Let $E$ be a real Banach space. A nonempty convex closed subset
$P\subset E$ is called a cone in $E$ if it satisfies:
\begin{enumerate}
  \item[$(i)$] {$P$ is closed, nonempty and $P\neq \{0\}$,}
  \item[$(ii)$] {$a,b\in \mathbb{R},$ $a,b\geq 0$ and $x,y \in P$ imply that $ax+by \in P,$}
  \item[$(iii)$] {$x \in P$ and $-x \in P$ imply that $x = 0.$}
\end{enumerate}
The space $E$ can be partially ordered by the cone $P\subset E$;
that is, $x \le y$ if and only if  $y-x \in P$. Also we write
$x\ll y$ if $y-x \in P^o$, where $P^o$ denotes the interior of $P$.\\
A cone $P$ is called normal if there exists a constant $K>0$ such
that $0\le x \le y$ implies $\|x\| \le K\|y\|$.\\
In the sequel we always suppose that $E$ is a real Banach
space, $P$ is a cone in $E$ with nonempty interior i.e. $P^o\neq \emptyset$ and $\leq$ is the partial ordering with
respect to $P$.
%---------------------------------------------------------------------------------------------------------------------
                    %definition 1.1
%---------------------------------------------------------------------------------------------------------------------
\begin{defn}
(\cite{Xian}) Let $X$ be a nonempty set.
Assume that the mapping $d:X\times X\rightarrow E$ satisfies
\begin{enumerate}
  \item[(i)] {$0\leq d(x,y)$ for all $x,y \in X$ and $d(x,y)=0 $ iff $x=y$}
  \item[(ii)] {$d(x,y)=d(y,x)$ for all $x,y \in X$}
  \item[(iii)] {$d(x,y)\leq d(x,z)+d(z,y)$ for all $x,y,z \in X$}.
\end{enumerate} Then $d$ is called a cone metric on $X$, and
$(X,d)$ is called a cone metric space.
\end{defn}
%---------------------------------------------------------------------------------------------------------------------
The study of fixed point theorems in such spaces is followed by some other mathematicians,
see \cite{AJ}-\cite{C2}. But the main question
was \textit{ "Are cone metric spaces a real generalization of metric spaces?."}
Throughout this paper we answer the question in the negative,  proving that
every cone metric space is metrizable  and the equivalent metric satisfies
the same contractive conditions as cone metric. So most of the fixed point
theorems which have been proved are the straightforward results from the metric case.
%---------------------------------------------------------------------------------------------------------------------
  \section{Main results}
%---------------------------------------------------------------------------------------------------------------------
                 % main theorem 2.1
%---------------------------------------------------------------------------------------------------------------------
\begin{thm}%\label{main}
For every cone metric $D:X\times X\rightarrow E$ there exists metric $d:X\times X\rightarrow \mathbb{R}^+$ which is equivalent to $D$ on $X$.
\end{thm}
{\bf\it Proof.} Define $d(x,y)=\inf\{\|u\|: D(x,y)\leq u\}.$
We shall to prove that $d$ is an equivalent metric to $D$.
If $d(x,y)=0$ then there exists $u_n$ such that $\|u_n\|\rightarrow 0 $ and $D(x,y)\leq u_n$. So $u_n\rightarrow 0 $ and consequently  for all $c\gg 0$ there exists $N\in \mathbb{N}$ such that $u_n\ll c$ for all $n\geq N.$ Thus for all $c\gg 0$, $0\leq D(x,y)\ll c$. Namely $x=y.$ \\
If $x=y$ then $D(x,y)=0$ which implies that $d(x,y)\leq \|u\|$ for all $0\leq u$. Put $u=0$ it implies $d(x,y)\leq \|0\|=0,$ on the other hand $0\leq d(x,y),$ therefore $d(x,y)=0.$ It is clear that $d(x,y)=d(y,x).$ To prove triangle inequality, for $x,y,z\in X$ we have,
$$\forall \varepsilon>0 \quad\exists u_1 \quad \|u_1\|<d(x,z)+\varepsilon,\quad D(x,z)\leq u_1,
$$
$$\forall \varepsilon>0 \quad\exists u_2 \quad \|u_2\|<d(z,y)+\varepsilon,\quad D(z,y)\leq u_2.
$$
But $D(x,y)\leq D(x,z)+D(z,y)\leq u_1+u_2,$ therefore
$$d(x,y)\leq \|u_1+u_2\|\leq \|u_1\|+\|u_2\|\leq d(x,z)+d(z,y)+2\varepsilon.$$
Since $\varepsilon>0$ was arbitrary so $d(x,y)\leq  d(x,z)+d(z,y).$\\
Now we shall prove that, for all $\{x_n\}\subseteq X$ and $x\in X$, $x_n\rightarrow x$ in $(X,d)$ if and only if $x_n\rightarrow x$ in $(X,D).$
We have
$$\forall n,m\in \mathbb{N} \quad \exists u_{nm}\quad such~ that\quad \|u_{nm}\|<d(x_n,x)+\frac{1}{m}, \quad D(x_n,x)\leq u_{nm}.
$$
Put $v_n:=u_{nn}$ then $\|v_n\|<d(x_n,x)+\frac{1}{n}$ and $D(x_n,x)\leq v_n.$
Now if $x_n\rightarrow x$ in $(X,d)$ then $d(x_n,x)\rightarrow 0$ and so $v_n\rightarrow 0$ too, therefore for all $c\gg 0$ there exists $N\in \mathbb{N}$ such that $v_n\ll c$ for all $n\geq N.$
This implies that $D(x_n,x)\ll c$ for all $n\geq N.$ Namely $x_n\rightarrow x$ in $(X,D)$.\\
Conversely, for every  real $\varepsilon>0 $, choose $c\in E$ with $c\gg 0$ and $\|c\|<\varepsilon$. Then there exists $N\in \mathbb{N}$ such that $D(x_n,x)\ll c$ for all $n\geq N.$ This means that for all $\varepsilon>0$ there exists $N\in \mathbb{N}$ such that $d(x_n,x)\leq \|c\|<\varepsilon$ for all $n\geq N.$
Therefore $d(x_n,x)\rightarrow 0$ as $n\rightarrow \infty$ so $x_n\rightarrow x$
in $(X,d)$.$Box$
\newline
%---------------------------------------------------------------------------------------------------------------------
                 %example 2.2
%---------------------------------------------------------------------------------------------------------------------
\begin{exa}
Let $0\neq a\in P\subseteq \mathbb{R}^n$ with $\|a\|=1$ and for every $x,y\in \mathbb{R}^n$ define
$$
D(x,y)=\left\{
                                                                \begin{array}{ll}
                                                                  a, & \hbox{$x\neq y$;} \\
                                                                  0, & \hbox{$x=y$.}
                                                                \end{array}
                                                              \right.
$$
Then $D$ is a cone metric on $\mathbb{R}^n$ and its equivalent metric $d$ is
$$d(x,y)=\left\{
                    \begin{array}{ll}
                      1, & \hbox{$x\neq y$;} \\
                      0, & \hbox{$x=y$,}
                    \end{array}
                  \right.
$$
which is discrete metric.
\end{exa}
%---------------------------------------------------------------------------------------------------------------------
                  %example 2.3
%---------------------------------------------------------------------------------------------------------------------
\begin{exa}
 Let $a,b\geq 0$ and consider the cone metric $D:\mathbb{R}\times \mathbb{R}\rightarrow \mathbb{R}^2 $ with $D(x,y)=(ad_1(x,y),bd_2(x,y))$ where $d_1,d_2$ are metrics on $\mathbb{R}$. Then its equivalent metric is
$d(x,y)=\sqrt{a^2+b^2}\|(d_1(x,y),d_2(x,y))\|.$ In particular if $d_1(x,y):=|x-y|$ and $d_2(x,y):=\alpha |x-y|$,
where $\alpha\geq 0$ then $D$ is the same famous cone metric which has been introduced in \cite[Example 1]{Xian}
and its equivalent metric is $d(x,y)=\sqrt{1+\alpha^2}|x-y|.$
\end{exa}
%---------------------------------------------------------------------------------------------------------------------
                 %example 2.4
%---------------------------------------------------------------------------------------------------------------------
\begin{exa}
For $q> 0,$ $ b>1$, $E=l^q$, $P=\{ \{x_n\}_{n\geq 1}:x_n\geq 0,\quad for~ all ~ n\}$ and $(X,\rho)$
a metric space, define $D:X\times X\rightarrow E$ which is the same cone metric as \cite[Example 1.3]{HSh} by
$$D(x,y)=\left\{ \left(\frac{\rho(x,y)}{b^n}\right)^{\frac{1}{q}}\right\}_{n\geq 1}.$$
Then its equivalent metric on $X$ is
$$d(x,y)=\left\|\left\{ \left(\frac{\rho(x,y)}{b^n}\right)^{\frac{1}{q}}\right\}_{n\geq 1} \right\|_{l^q}=\left(\sum_{n=1}^{\infty}\frac{\rho(x,y)}{b^n}\right)^{\frac{1}{q}}
=\left(\frac{\rho(x,y)}{b-1}\right)^{\frac{1}{q}}.
$$
\end{exa}
%---------------------------------------------------------------------------------------------------------------------
                 % Lemma 2.5+
%---------------------------------------------------------------------------------------------------------------------
\begin{lem}\label{b}
Let $D,D^*:X\times X\rightarrow \Bbb E$ be cone metrics, $d,d^*:X\times X\rightarrow \mathbb{R}^+$ their equivalent metrics respectively and $T:X\rightarrow X$ a self map. If  $D(Tx,Ty)\leq D^*(x,y)$, then $d(Tx,Ty)\leq d^*(x,y).$
\end{lem}
{\bf\it Proof.} By the definition of $d^*$,
$$\forall \varepsilon>0 \quad\exists v \quad \|v\|<d^*(x,y)+\varepsilon,\quad D^*(x,y)\leq v.
$$
Therefore if $D(Tx,Ty)\leq D^*(x,y)\leq v,$
then we have $$d(Tx,Ty)\leq \| v\|\leq d^*(x,y)+\varepsilon.$$
Since $\varepsilon>0$ was arbitrary so $d(Tx,Ty)\leq d^*(x,y).$$\Box$
%\newline
%---------------------------------------------------------------------------------------------------------------------
                 %example 2.6+
%---------------------------------------------------------------------------------------------------------------------
\begin{exa}
Let $E:=\Bbb R$, $P:=\Bbb R^+$ and $D:X\times X\rightarrow E$ be a cone metric, $d:X\times X\rightarrow \mathbb{R}^+$ its equivalent metric, $T:X\rightarrow X$ a self map and $\varphi:\Bbb R^+\rightarrow \Bbb R^+$ defined by $\varphi (x)=\frac{x}{1+x}$. If  $D^*:=\varphi(D)$, then its equivalent metric is $d^*=\varphi(d)$, and if, $D(Tx,Ty)\leq \varphi (D(x,y))=\frac{D(x,y)}{1+D(x,y)},$
then by lemma \ref{b}, $d(Tx,Ty)\leq \varphi (d(x,y))=\frac{d(x,y)}{1+d(x,y)}.$
\end{exa}
%---------------------------------------------------------------------------------------------------------------------
                %definition 2.7
%---------------------------------------------------------------------------------------------------------------------
\begin{defn}
A self map $\varphi$ on normed space $X$ is bounded if
$$\|\varphi\|:=\sup_{0\neq x\in X}\frac{\|\varphi(x)\|}{\|x\|}<\infty.
$$
\end{defn}
%---------------------------------------------------------------------------------------------------------------------
                 % theorem 2.8
%---------------------------------------------------------------------------------------------------------------------
\begin{thm}
Let $D:X\times X\rightarrow  E$ be a cone metric, $d:X\times X\rightarrow \mathbb{R}^+$ its equivalent metric, $T:X\rightarrow X$ a self map and $\varphi:P\rightarrow P$ a bounded map, then there exists
$\psi:\mathbb{R}^+\rightarrow \mathbb{R}^+$
such that $D(Tx,Ty)\leq \varphi(D(x,y))$ for every $x,y\in X$ implies $d(Tx,Ty)\leq \psi(\|D(x,y)\|)$ for all $x,y\in X$. Moreover if $\psi$ is decreasing map or $\varphi$ is linear and increasing map then,
$d(Tx,Ty)\leq \psi(d(x,y))$ for all $x,y\in X$.
\end{thm}
{\bf\it Proof.} Put $\psi(t):=\sup_{0\neq x\in P}\left\|\varphi\left(\frac{t}{\|x\|}x\right )\right\|$ for all $t\in \Bbb R^+.$
So $\|\varphi(x)\|\leq\psi(\|x\|)$  for all $x\in P$.
Therefore if $D(Tx,Ty)\leq \varphi(D(x,y)),$
then we have $d(Tx,Ty)\leq \| \varphi(D(x,y))\|\leq \psi(\|D(x,y)\|).$
Now if $\psi$ be decreasing map, by the definition of $d$ we have $d(x,y)\leq \|D(x,y)\|$, so
$$d(Tx,Ty)\leq  \psi(\|D(x,y)\|)\leq \psi(d(x,y)).$$  If $\varphi$ be a
linear increasing map then $\psi(t)=t\|\varphi\|.$ The definition of $d$ implies
$$\forall \varepsilon>0 \quad\exists v \quad \|v\|<d(x,y)+\varepsilon,\quad D(x,y)\leq v.
$$
Therefore if $D(Tx,Ty)\leq \varphi(D(x,y))\leq \varphi(v),$
then we have
$$d(Tx,Ty)\leq \| \varphi(v)\|\leq \psi(\|v\|)\leq \psi(d(x,y))+\psi(\varepsilon).$$
Since $\varepsilon>0$ was arbitrary and $\psi(\varepsilon)\rightarrow 0$ as $\varepsilon\rightarrow 0$,  so $d(Tx,Ty)\leq  \psi(d(x,y)).$$\Box$

In the following summary of our results are listed.

%====================================================
\begin{cor}
Let $D,D^*$ be cone metrics, $d.d^*$ their equivalent metrics, $T :X\rightarrow X$ a map, $\lambda\in[0,\frac{1}{2})$ and $\alpha,\beta\in [0,1)$. For $x, y \in X$,
\begin{enumerate}
\item[1] $D(Tx, Ty)\leq \alpha D(x,y) \Rightarrow d(Tx, Ty)\leq \alpha d(x, y).$
  \item[2] $D(Tx, Ty)\leq \lambda (D(Tx, x) + D(Ty, y)) \Rightarrow d(Tx, Ty)\leq \lambda (d(Tx, x) + d(Ty, y)).$
  \item[3] $D(Tx, Ty)\leq \lambda (D(Tx, y) + D(Ty, x)) \Rightarrow d(Tx, Ty)\leq \lambda (d(Tx, y) + d(Ty, x)).$
  \item[4] $D(Tx, Ty)\leq \alpha D(x, y) + \beta D(Tx, y) \Rightarrow d(Tx, Ty)\leq \alpha d(x, y) +\beta d(Tx, y).$
      \item[5] There exists $u\in \{D(x, y); D(x, Tx); D(y, Ty); \frac{1}{2}[
D(x, Ty)]+ D(y, Tx)]\}$ such that $D(Tx, Ty)\leq \alpha u$ where $\alpha \in(0, 1),$
 then $$d(Tx, Ty)\leq \alpha\max\{d(x, y); d(x, Tx); d(y, Ty); \frac{1}{2}[d(x, Ty)]+ d(y, Tx)]\}.$$
  \item[6] There exists $u\in \{D(x, y); D(x, Tx); D(y, Ty); \frac{1}{2}
D(x, Ty);\frac{1}{2}D(y, Tx)\}$ such that $D(Tx, Ty)\leq \beta u$ where $\beta \in(0, 1),$ then $$d(Tx, Ty)\leq \beta\max\{d(x, y); d(x, Tx); d(y, Ty); \frac{1}{2}d(x, Ty);\frac{1}{2}d(y, Tx)\}.$$
  \item[7] There exists $u\in \{D(x, y); \frac{1}{2}[D(x, Tx)+ D(y, Ty)]; \frac{1}{2}[D(x, Ty)+ D(y, Tx)]\}$ such that $D(Tx, Ty)\leq \beta u$ where $\beta \in(0, 1),$ then $$d(Tx, Ty)\leq \beta\max\{d(x, y); \frac{1}{2}[d(x, Tx)+ d(y, Ty)]; \frac{1}{2}[d(x, Ty)+ d(y, Tx)]\}.$$
  \item[8]If $$D(Tx, Ty)\leq a_1D(x, y)+a_2 DSx, Tx)+a_3 D(y, Ty)+a_4D(x, Ty)+a_5D(y, Tx),$$ then $$d(Tx, Ty)\leq a_1d(x, y)+a_2 d(x, Tx)+a_3 d(y, Ty)+a_4d(x, Ty)+a_5d(y, Tx)$$ where $\sum_{i=1}^5a_i< 1.$
  \item[9] If there exists $$u\in \{D(x, y); D(x, Tx); D(y, Ty); D(x, Ty);D(y, Tx)\}$$ such that  $D(Tx, Ty)\leq \frac{\beta}{2}u,$ then $$d(Tx, Ty)\leq \frac{\beta}{2}\max\{d(x, y); d(x, Tx); d(y, Ty); d(x, Ty);d(y, Tx)\}$$ where $\beta \in(0, 1).$
  \item[10] If $$D(Tx, Ty)\leq a_1D(x, y)+a_2 D(x, Tx)+a_3 D(y, Ty)+a_4[D(x, Ty)+D(y, Tx)],$$ then $$d(Tx, Ty)\leq a_1d(x, y)+a_2 d(x, Tx)+a_3 d(y, Ty)+a_4[d(x, Ty)+d(y, Tx)]$$ where $a+1+a_2+a_3+2a_4< 1.$
      \item[11] There exist $m,n\in\mathbb{N}$ and $ k\in[0,1)$ such that $$D(T^mx,T^ny)\leq kD(z,t)$$ for all $x,y\in X$, $z\neq t$ and $z,t\in \{x,y,T^px,T^qy\}$ where $1\leq p\leq m$ and $1\leq q\leq n,$ then $$d(T^mx,T^ny)\leq k d(z,t).$$
          \item[12] If $D(Tx,Ty)\leq D^*(x,y),$ then $d(Tx,Ty)\leq d^*(x,y).$
      \end{enumerate}
\end{cor}
%---------------------------------------------------------------------------------------------------------------------
                   %Acknowledgments
%---------------------------------------------------------------------------------------------------------------------
{\bf Acknowledgements}\newline
This paper has been supported by the Islamic Azad University, Zanjan Branch, Zanjan, Iran. The first author would like to thank this support. And the authors would  like to thank {\it Professor Brailey Sims} for his helpful advise which led them to  present  this paper.

%======================================================================================

\end{document}